\documentclass[10pt]{article}
\usepackage{latexsym}
\usepackage{amsbsy}
\usepackage{amssymb}
\usepackage{amsfonts}
\usepackage{amsmath}
\usepackage[latin1]{inputenc}

 \newtheorem{thm}{Theorem}[section]
 \newtheorem{defin}[thm]{Definition}
 \newtheorem{lem}[thm]{Lemma}
 \newtheorem{prop}[thm]{Proposition}
 \newtheorem{cor}[thm]{Corollary}
 \newtheorem{rem}[thm]{Remark}
 \newtheorem{ex}[thm]{Example}

 \newcommand{\bthm}{\begin{thm}}
 \newcommand{\ethm}{\end{thm}}
 \newcommand{\bd}{\begin{defin}}
 \newcommand{\ed}{\end{defin}}
 \newcommand{\blem}{\begin{lem}}
 \newcommand{\elem}{\end{lem}}
 \newcommand{\bcor}{\begin{cor}}
 \newcommand{\ecor}{\end{cor}}
 \newcommand{\bprop}{\begin{prop}}
 \newcommand{\eprop}{\end{prop}}
 \newcommand{\brem}{\begin{rem} \rm}
 \newcommand{\erem}{\end{rem}}
 \newcommand{\bex}{\begin{ex} \rm}
 \newcommand{\eex}{\end{ex}}

 \newcommand{\pr}{\noindent{\bf Proof. }}
 \newcommand{\ep}{\nolinebreak{\hspace*{\fill}$\Box$ \vspace*{0.25cm}}}

 \newcommand{\beq}{\begin{equation}}
 \newcommand{\eeq}{\end{equation} }

 \newcommand{\bea}{\begin{eqnarray}}
 \newcommand{\eea}{\end{eqnarray}}

 \newcommand{\beas}{\begin{eqnarray*}}
 \newcommand{\eeas}{\end{eqnarray*}}

 \newcommand{\beqs}{\begin{equation*}}
 \newcommand{\eeqs}{\end{equation*}}

 \newcommand{\bi}{\begin{itemize}}
 \newcommand{\ei}{\end{itemize}}

 \newcommand{\ben}{\begin{enumerate}}
 \newcommand{\een}{\end{enumerate}}

 \newcommand{\ba}{\begin{array}}
 \newcommand{\ea}{\end{array}}

 \newcommand{\R}{\mathbb R}

 \newcommand{\cC}{\ensuremath{{\cal C}}}

 \newcommand{\cL}{\ensuremath{{\cal L}}}

 \newcommand{\cU}{\ensuremath{{\cal U}}}

 \newcommand{\pd}{\partial}
 \renewcommand{\d}{\ensuremath{\partial}}

 \newcommand{\eps}{\varepsilon}

 \newcommand{\al}{\alpha}
 
\begin{document}

 \title{Generalized Hamilton's Principle with Fractional Derivatives}
 \author{Teodor Atanackovi\'c
         \footnote{Faculty of Technical Sciences, Institute of Mechanics, University of Novi Sad,
         Trg Dositeja Obradovi\'ca 6, 21000 Novi Sad, Serbia,
         Electronic mail: atanackovic@uns.ac.rs}\\
         Sanja Konjik
         \footnote{Faculty of Sciences, Department of Mathematics and Informatics, University of Novi Sad,
         Trg Dositeja Obradovi\'ca 4, 21000 Novi Sad, Serbia,
         Electronic mail: sanja.konjik@dmi.uns.ac.rs}\\
         Ljubica Oparnica
         \footnote{Institute of Mathematics, SANU,
         Kneza Mihaila 36, 11000 Belgrade, Serbia,
         Electronic mail: ljubicans@sbb.rs}\\
         Stevan Pilipovi\'c
         \footnote{Faculty of Science, Department of Mathematics and Informatics, University of Novi Sad,
         Trg Dositeja Obradovi\'ca 4, 21000 Novi Sad, Serbia,
         Electronic mail: pilipovic@dmi.uns.ac.rs}\\
       }

 \date{}
 \maketitle

 \begin{abstract}
 We generalize Hamilton's principle with fractional derivatives
 in Lagrangian $L(t,y(t),{}_0D_t^\al y(t),\alpha)$
 so that the function $y$ and the order of fractional
 derivative $\alpha$ are varied in the minimization procedure.
 We derive stationarity conditions and discuss
 them through several examples.

 \vskip5pt
 \noindent
 {\bf Mathematics Subject Classification (2000):}
 Primary: 49K05; Secondary: 26A33, 70H25
 \vskip5pt
 \noindent
 {\bf Keywords:} Hamilton's principe, fractional derivatives,
 optimality conditions, stationary points
 \end{abstract}

\section{Introduction}

 Hamilton's principle is one of the basic principles of Physics.
 Anthony \cite{Anthony} states: "In theoretical physics, a theory
 is often considered to be complete if its variational principle
 in the sense of Hamilton is known".
 When a Hamilton's principle is known, the whole information
 concerning the processes of a particular system is included into
 its Lagrangian.
 When known, the Hamilton's principe could be used in many
 different ways. For example, it could serve as the basis for
 obtaining first integrals via N\"{o}ther's theory, or to generate
 approximate solutions to the relevant system of equations by the
 use of Ritz procedure (cf.\ \cite{VujanovicAtanackovic}).

 In this paper we investigate necessary conditions
 for solutions of fractional variational problems (Euler-Lagrange
 equations). Such investigations has been initiated in
 \cite{Riewe96, Riewe97},
 and continued in \cite{Agrawal02} - \cite{AKP-fracEL} and
 \cite{BalTru08}; see also \cite{Rekh} for
 the importance of introducing fractional derivatives into
 the Lagrangian density of Hamilton's principle.
 In general, we refer to
 \cite{Dacorogna} - \cite{Jumarie07},
 \cite{MillRoss, OldSpa, Podlubny},
 \cite{SamkoKM} - \cite{WestBolognaGrigolini},
 for different aspects of the calculus of variations and fractional
 calculus, motivations and applications.

 When a fractional variational problem is studied,
 a natural question arises how one can choose $\alpha$, the order
 of fractional derivative, in order to
 achieve the minimal value of a functional under consideration.
 More generally, one can address the same question for any problem
 involving fractional operators.
 Usually, in application, a good choice of $\alpha$
 is obtained by experiments, numerical methods or
 computational simulations. However, experimental results give
 different values for $\alpha$ within certain interval.
 In this paper we
 propose a method which a priori gives values for $\alpha$
 which optimize the considered variational problem following the
 fundamental minimization principle of Hamilton's action.
 In fact, we address the question of finding stationary
 points for Hamilton's action integral with fractional Lagrangian
 in a more general setting. Namely, we allow the stationarity
 of the action integral with respect to a set of admissible
 functions \textit{and} with respect to the order of fractional
 derivatives, appearing in the Lagrangian.
 Up to our knowledge the problem when both $y$
 and $\alpha$ are varied has
 not been analyzed so far. It leads to stationarity conditions
 as a basis for generalized Hamilton's principle
 for the action integral
 \begin{equation} \label{eq:fvp-alpha}
 I[y,\alpha]: = \int_{0}^{b}L( t,y( t),{}_{0}D_{t}^{\alpha }y,\alpha)\,
 dt,\quad y\in \cU,\alpha\in A:=[0,\alpha_0],\alpha_0\leq 1,
 \end{equation}
 where $\cU$ is a set of admissible functions: Find
 \begin{equation}\label{eq:min}
 \min_{(y,\alpha)\in \cU\times A} I[y, \alpha]
 \qquad \mbox{or}
 \end{equation}
 \begin{equation}\label{eq:minmin}
 \min_{\alpha\in A}(\min_{y\in\cU} I[y, \alpha])
 \qquad \mbox{or}
 \end{equation}
 \begin{equation}\label{eq:minimin}
 \min_{y\in\cU}( \min_{\alpha\in A} I[y, \alpha]).
 \end{equation}
 In this paper we analyze stationarity conditions for
 (\ref{eq:min}) and (\ref{eq:minmin}).
 Stationarity conditions with respect to $\alpha$
 in (\ref{eq:minimin}) are more difficult,
 and, contrary to (\ref{eq:minmin}), this case is less natural in
 applications. Note that in (\ref{eq:min}), (\ref{eq:minmin}) and
 (\ref{eq:minimin}) one can look for maximums instead of minimums.
 So, our general problem is determination of stationary points.

 So far, parameter $\alpha$, the order of fractional derivative,
 has been determined experimentally (cf.\ e.g. \cite{Rogers83}).
 This approach offers a rational way for choosing the precise
 $\alpha$.

 The paper is organized as follows. To the end of Introduction, we
 recall definitions and properties of fractional derivatives.
 In Section \ref{sec:Hamilton} we present a framework in which we
 shall study variational problems (\ref{eq:min}) and
 (\ref{eq:minmin}). Then in Section \ref{SecOptCond}
 we derive stationarity conditions for (\ref{eq:min}).
 Section \ref{SecEqProbl} is devoted to additional
 assumptions which provide equivalence of problems
 (\ref{eq:min}) and (\ref{eq:minmin}).
 Results which are obtained in previous sections are illustrated
 by several examples in Section \ref{SecEx}.
 Moreover, examples of this section give further motivation for
 our investigation.
 In the last remark of
 Section \ref{SecEx} we propose a new formulation of a
 fractional variational problem.

 Throughout this paper we shall use the following notation. The
 mapping $(t,\alpha)\mapsto {}_0D_t^\alpha(y)$, which defines
 the left Riemann-Liouville fractional derivative of
 order $\alpha$, will be denoted by ${}_0D_t^\alpha y$, or
 by ${}_0D_t^\alpha y(t)$. Recall,
 $$
 {}_{0}D_{t}^{\alpha}y:=\frac{1}{\Gamma(1-\alpha)}\frac{d}{dt}
 \int_{0}^{t}\frac{y(\tau)}{(t-\tau)^{\alpha}} \, d\tau,
 \quad t\in[0,b],\, 0\leq\alpha <1,
 $$
 where $\Gamma$ is the Euler gamma function,
 and its existence is provided whenever
 \begin{equation} \label{eq:acfd}
 [0,b]\ni t\mapsto
 \int_0^t \frac{y(\tau)}{(t-\tau)^{\alpha}}\, d\tau
 \end{equation}
 is an absolutely continuous function. Recall, the space of
 absolutely continuous functions is denoted by
 $AC([0,b])$ and it is supplied with the norm
 $||f||=\sup_{x\in[0,b]} |f(x)|$
 (clearly, it is not a Banach space).
 For example, (\ref{eq:acfd}) is absolutely continuous if
 $y\in AC([0,b])$.
 However, there are some cases when with less regularity in $y$ we
 still have a well-defined operator of fractional differentiation (cf.\
 \cite{SamkoKM}). For instance, ${}_{0}D_{t}^{\alpha}y$ exists for
 functions with integrable singularities
 (a continuous and locally integrable function $f$ in $(0,b]$ has
 an integrable singularity
 at $\tau=0$ of order $r<1$ if
 $\lim_{\tau\to 0}\tau^r f(\tau)\not= 0$).
 In particular we can take $y(t)=t^{-\mu}$,
 $t\in(0,b]$ (for any $b>0$), $0<\mu <1$.
 Then we obtain the so-called Euler formula
 (cf.\ \cite[(2.26)]{SamkoKM})
 $$
 {}_0D_t^\alpha t^{-\mu}
 =\frac{\Gamma(1-\mu)}{\Gamma(1-\mu-\alpha)}
 \frac{1}{t^{\mu+\alpha}}, \quad t\in(0,b].
 $$

 The right Riemann-Liouville fractional derivative of order
 $\alpha$ is defined as
 $$
 {}_{t}D_{b}^{\alpha}y:=\frac{1}{\Gamma(1-\alpha)}
 \left(-\frac{d}{dt}\right)
 \int_{t}^{b}\frac{y(\tau)}{(\tau-t)^{\alpha}} \, d\tau,
 \quad t\in[0,b],\, 0\leq\alpha <1.
 $$
 The conditions for its existence are similar as in the case of
 the left fractional derivative.

 In the sequel we shall consider cases involving both fractional
 derivatives and work with integrable functions for which these
 derivatives (or one of them) exist. In such cases notation
 ${}_{0}D_{t}^{\alpha}y$, resp.\ ${}_{t}D_{b}^{\alpha}y$, $t\in[0,b]$,
 means that $y$ and ${}_{0}D_{t}^{\alpha}y$, resp.\
 ${}_{t}D_{b}^{\alpha}y$, are considered as integrable functions
 which can take values $+\infty$ or $-\infty$ at some points.

 Recall (\cite{Nakhushev}), that
 ${}_0D_t^\alpha y\to y'$ and ${}_tD_b^\alpha y\to -y'$ in $\cC([0,b])$,
 as $\alpha\to 1^{-}$, whenever $y\in \cC^1([0,b])$.

 Also, we shall make use of Caputo fractional derivatives.
 The left, resp.\ right, Caputo fractional derivative
 of order $\alpha\in[0,1)$ is defined as
 $$
 {}_{0}^{C}D_{t}^{\alpha}y: = \frac{1}{\Gamma(1-\alpha)}
 \int_{0}^{t} \frac{y'(\tau)}{(t-\tau)^{\alpha}}\, d\tau
 \quad \mbox{resp.} \quad
 {}_{t}^{C}D_{b}^{\alpha}y: = \frac{1}{\Gamma(1-\alpha)}
 \int_{t}^{b} \frac{-y'(\tau)}{(\tau-t)^{\alpha}}\, d\tau.
 $$
 One can show that for $y\in AC([0,b])$ and $t\in [0,b]$,
 $$
 {}_0D_t^\alpha y= {}_{0}^{C}D_{t}^{\alpha}y
 + \frac{1}{\Gamma (1-\alpha)}
 \frac{y(0)}{(t-a)^\alpha},
 \quad
 {}_tD_b^\alpha y= {}_{t}^{C}D_{b}^{\alpha}y
 + \frac{1}{\Gamma (1-\alpha)}
 \frac{y(b)}{(b-t)^\alpha},
 $$
 (cf.\ e.g. \cite{KilSriTru}).
 Therefore, ${}_{0}D_{t}^{\alpha}y={}_{0}^{C}D_{t}^{\alpha}y$,
 resp.\ ${}_{t}D_{b}^{\alpha}y={}_{t}^{C}D_{b}^{\alpha}y$,
 whenever $y(0)=0$, resp.\ $y(b)=0$.

 \section{Formulation of the problem}
 \label{sec:Hamilton}

 We investigate stationary points of (\ref{eq:fvp-alpha})
 for $\alpha\in[0,\alpha_0]$ and
 all admissible functions $y$, whose properties will be specified
 in the sequel.
 We shall distinguish two cases: $\alpha_0$ strictly less than $1$
 and $\alpha_0=1$.
 In the case $\alpha_0<1$, set
 $$
 \mathcal{U}_l:=\{y\in L^1([0,b])\,|\,
 {}_0D_t^\alpha y\in L^1([0,b])\}.
 $$
 Obviously, $AC([0,b])$ is a subset of $\cU_l$.
 In the case $\alpha_0=1$ we assume that
 $y\in\mathcal{U}_l$ and that, in addition,
 ${}_0D_t^1 y$ exists, and ${}_0D_t^1 y=y'$
 is an integrable function.
 Let us note that one can consider $\mathcal{U}_l$ defined with
 $L^p([0,b])$ (or their subspaces) instead of $L^1([0,b])$
 (see Remark \ref{rem:lin ops}).

 In general, we shall use notation
 \beq \label{admisset}
 \mathcal{U}:=\{y\in \cU_l\,|\,
 y \mbox{ satisfies specified boundary conditions}\}.
 \eeq
 We shall sometimes write $\mathcal{U}$ also for $\mathcal{U}_l$
 (then the set of specified boundary conditions is empty).

 In the sequel, Lagrangian $L(t,y(t),{}_{0}D_{t}^{\alpha}y,\alpha)$
 (Lagrangian density, in Physics) satisfies:
 \beq \label{eq:conditions on L}
 \left.\ba{c}
 L\in\cC^1([0,b]\times \R\times \R\times [0,1])\\
 \mbox{and}\\
 t\mapsto\pd_3 L(t,y(t),{}_0D_t^\alpha y,\alpha)
 \in \cU_r, \mbox{ for every } y\in \cU_l
 \ea\right\}
 \eeq
 where $\mathcal{U}_r:=\{y\in L^1([0,b])\,|\,
 {}_tD_b^\alpha y\in L^1([0,b])\}$.

 Recall, our generalization of Hamilton's principle is realized
 through the determination of
 $(y^{\ast},\alpha^{\ast})\in\mathcal{U}\times A$ such that
 \begin{equation} \label{eq:fvp-alpha-min}
 \underset{(y,\alpha) \in \mathcal{U}\times A}{\min}
 \int_{0}^{b} L(t,y(t),{}_{0}D_{t}^{\alpha}y,\alpha)\, dt
 =
 \int_{0}^{b} L( t,y^{\ast}(t),{}_{0}D_{t}^{\alpha^{\ast}}y^{\ast},
 \alpha^{\ast})\,dt.
 \end{equation}

 There are two special cases of (\ref{eq:fvp-alpha-min}).
 The first one is obtained when $A=\{1\}$.
 Then, since $\left.{}_0D_t^\alpha y\right\vert_{\alpha=1}
 =y'(t)$ for $y\in\cC^1([0,b])$,
 the solution $y^{\ast}$ of (\ref{eq:fvp-alpha-min})
 satisfies the classical Euler-Lagrange equation
 $$
 \frac{d}{dt}\frac{\partial L}{\partial y'}-
 \frac{\partial L}{\partial y} =0.
 $$

 If $A$ has a single element $A=\{\alpha\}$, $0< \alpha < 1$, then
 $\underset{(y,\alpha)\in\mathcal{U}\times \{\alpha\}}{\min}
 I[y,\alpha]$ leads to the fractional Euler-Lagrange equation
 (cf.\ (\cite{Agrawal02, AKP-fracEL}))
 $$
 {}_{t}D_{b}^{\alpha}\frac{\partial L}{\partial {}_{0}D_{t}^{\alpha}y}
 +\frac{\partial L}{\partial y}=0.
 $$

 We proceed with finding stationary points related to
 (\ref{eq:fvp-alpha}).

 \section{Optimality conditions}
 \label{SecOptCond}

 A necessary condition for the existence of solutions to variational
 problem (\ref{eq:fvp-alpha-min}) is given in the following
 proposition.

 \bprop \label{prop:ELeqs}
 Let $L$ satisfy (\ref{eq:conditions on L}).
 Then a necessary condition that functional
 (\ref{eq:fvp-alpha}) has an extremal point at
 $(y^{\ast},\alpha^{\ast})\in\cU\times A$
 is that $(y^{\ast},\alpha^{\ast})$ is
 a solution of the system of equations
 \begin{eqnarray}
 \frac{\partial L}{\partial y}+{}_{t}D_{b}^{\alpha}
 \frac{\partial L}{\partial{}_{0}D_{t}^{\alpha}y}
 &=& 0,  \label{eq:ELeqs-y} \\
 \int_{0}^{b} \left(\frac{\partial L}{\partial{}_{0}D_{t}^{\alpha }y}G( y,\alpha)
 +\frac{\partial L}{\partial \alpha }\right)\, dt
 &=& 0,  \label{eq:ELeqs-alpha}
 \end{eqnarray}
 where
 $$
 G(y,\alpha)=\frac{\partial{}_{0}D_{t}^{\alpha}y}{\partial\alpha}
 =\frac{d}{dt}(f_{1}\ast_t y)(t,\alpha), \,
 f_{1}(t,\alpha)=\frac{1}{t^{\alpha}\Gamma(1-\alpha)}
 [\psi(1-\alpha)-\ln t], \, t>0,
 $$
 with the Euler function $\psi(z)=\frac{d}{dz}\ln\Gamma(z)$,
 and
 $(f_{1}\ast_t y)(t,\alpha)
 =\int_{0}^{t}f_{1}(\tau,\alpha)y(t-\tau)\, d\tau$.
 \eprop

 \pr
 Let $(y^{\ast},\alpha^{\ast})$ be an element of
 $\mathcal{U}\times A$
 for which $I[y,\alpha]$ has an extremal value. Let
 $y(t) =y^{\ast}(t) +\varepsilon_{1}f(t)$,
 $\alpha =\alpha^{\ast}+\varepsilon_{2}$,
 $\varepsilon_1,\varepsilon_2\in\R$,
 with $f\in\cU_l$, and the boundary conditions
 on $f$ are specified so that the varied path
 $y^{\ast}+\varepsilon_{1}f$ is an element of $\mathcal{U}$.
 Then $I[y,\alpha]=I[y^{\ast}+\varepsilon_{1}f,
 \alpha^{\ast}+\varepsilon_{2}]=:I(\varepsilon_{1},
 \varepsilon_{2})$.
 A necessary condition for an extremal value of $I[y,\alpha]$ is
 $$
 \left.\frac{\partial I(\varepsilon_{1},\varepsilon_{2})}{\partial\varepsilon_{1}}
 \right\vert_{\varepsilon_{1}=0,\varepsilon_{2}=0}=0,
 \qquad
 \left.\frac{\partial I(\varepsilon_{1},\varepsilon_{2})}{\partial\varepsilon_{2}}
 \right\vert_{\varepsilon_{1}=0,\varepsilon_{2}=0}=0.
 $$
 Therefore we obtain
 \begin{eqnarray}
 \int_{a}^{b}\left(\frac{\partial L}{\partial y}f(t)
 +\frac{\partial L}{\partial{}_{0}D_{t}^{\alpha}y}\,
 {}_{0}D_{t}^{\alpha }f\right)\, dt &=&0, \label{eq:11-1} \\
 \int_{0}^{b}\left(\frac{\partial L}{\partial {}_{0}D_{t}^{\alpha}y}\,
 \frac{\partial {}_{0}D_{t}^{\alpha}y}{\partial\alpha}
 +\frac{\partial L}{\partial \alpha }\right)\, dt &=&0.  \label{eq:11-2}
 \end{eqnarray}
 Applying the fractional integration by parts formula
 (cf.\ \cite{KilSriTru}):
 \begin{equation} \label{eq:frac int by parts}
 \int_{0}^{b} g(t)\cdot {}_{0}D_{t}^{\alpha }f(t)\, dt
 = \int_{0}^{b}\, f(t)\cdot {}_{t}D_{b}^{\alpha }g(t)\, dt,
 \end{equation}
 to (\ref{eq:11-1}), (\ref{eq:11-1}) is transformed
 to
 $$
 \int_{0}^{b} \left(\frac{\partial L}{\partial y}
 +{}_{t}D_{b}^{\alpha }\frac{\partial L}{\partial{}_{0}D_{t}^{\alpha}y}
 \right)f(t)\, dt=0.
 $$
 From this equation, using the fundamental lemma of
 the calculus of variations
 (see \cite[p.\ 115]{Dacorogna}), we conclude that
 condition (\ref{eq:ELeqs-y}) holds for the
 optimal values $y^{\ast}$ and $\alpha^{\ast}$.
 In (\ref{eq:ELeqs-alpha}) the term
 $\frac{\partial {}_{0}D_{t}^{\alpha}y}{\partial\alpha}$
 is transformed by the use of expression
 \begin{eqnarray}
 \frac{\partial {}_{0}D_{t}^{\alpha}y}{\partial\alpha}
 &=&\psi(1-\alpha)\, {}_{0}D_{t}^{\alpha}y -
 \frac{1}{\Gamma(1-\alpha)}\frac{d}{dt}
 \int_{0}^{t}
 \frac{\ln(t-\tau) y(\tau)}{(t-\tau)^{\alpha}}\, d\tau \notag \\
 &=&\frac{d}{dt}(f_{1}\ast_t y) (t,\alpha) \notag \\
 &=& G(y,\alpha), \quad (y,\alpha)\in\cU\times A \label{eq:14}
 \end{eqnarray}
 (cf.\ \cite[p.\ 592]{AOP07}).
 We obtain (\ref{eq:ELeqs-alpha}) by substituting (\ref{eq:14})
 into (\ref{eq:11-2}).
 \ep

 \brem
 In general, in  solving  equations
 (\ref{eq:ELeqs-y}) and (\ref{eq:ELeqs-alpha}), the most delicate
 task is the calculation of expression
 $\frac{\partial{}_{0}D_{t}^{\alpha}y}{\partial \alpha}$.
 Although its general form (\ref{eq:14}) has been derived
 in \cite[p.\ 592]{AOP07}, various difficulties can still appear.
 We illustrate this by examples in Section \ref{SecEx}.
 However, the simplified form of
 $\frac{\partial{}_{0}D_{t}^{\alpha}y}{\partial \alpha}$,
 in some special cases, is important.

 Already in
 \cite{AOP07} it has been shown that for $y\in AC([0,b])$
 \begin{eqnarray}
 \left.\frac{\partial {}_{0}D_{t}^{\alpha}y}{\partial\alpha}
 \right\vert_{\alpha =0^{+}}
 &=& -(\gamma +\ln t) y(0) -\int_{0}^{t}(\gamma +\ln \tau)
 y( t-\tau)\,d\tau \notag \\
 &=& -(\gamma +\ln t) y(t)
 + \int_{0}^{t}\frac{y(t)-y(t-\tau)}{\tau}\, d\tau, \label{16bbb}
 \end{eqnarray}
 where $\gamma =0.5772156...$ is the Euler constant. (Another form
 of $\left.\frac{\partial {}_{0}D_{t}^{\alpha}y}{\partial\alpha}
 \right\vert_{\alpha =0^{+}}$ is also given in
 \cite[p.\ 111]{WestBolognaGrigolini}.)

 Let us obtain a simplified form
 of $\frac{\partial {}_{0}D_{t}^{\alpha}y}{\partial\alpha}$ at
 $\alpha=1^{-}$. In order to do that we use the method
 proposed in \cite{TarasovZaslavsky06}.
 We recall the expansion of $(t-\tau)^{\eps}/\Gamma(1+\eps)$ with
 respect to $\eps$, at $\eps=0$, with $\tau<t$ (cf.\
 \cite[p.\ 401]{TarasovZaslavsky06}), which will be used in the
 sequel:
 \begin{equation} \label{eq:expansion}
 \frac{(t-\tau)^{\eps}}{\Gamma(1+\eps)}
 = \frac{e^{\eps\ln(t-\tau)}}{\Gamma(1+\eps)}
 = 1+\eps(\gamma+\ln(t-\tau))+ o(\eps).
 \end{equation}
 Assume now that $y\in C^{2}([0,b])$.
 Then, as in \cite{Podlubny}, for $t\in[0,b]$,
 \begin{eqnarray*}
 {}_{0}D_{t}^{\alpha }y &=& \frac{1}{\Gamma(1-\alpha)}
 \frac{d}{dt} \int_{0}^{t}\frac{y(\tau)}{(t-\tau)^{\alpha }}\, d\tau \\
 &=& \frac{y(0)}{\Gamma(1-\alpha) t^{\alpha}}
 +\frac{1}{\Gamma(1-\alpha)}
 \int_{0}^{t}\frac{y^{(1)}(\tau)}{(t-\tau)^{\alpha}}\, d\tau \\
 &=& \frac{y(0)}{\Gamma(1-\alpha) t^{\alpha }}+
 \frac{y^{(1)}(0)}{\Gamma(2-\alpha)t^{\alpha -1}}
 +\frac{1}{\Gamma(2-\alpha)}\int_{0}^{t}
 (t-\tau)^{1-\alpha}y^{(2)}(\tau)\, d\tau.
 \end{eqnarray*}
 With $\alpha=1-\varepsilon$ we obtain
 \begin{equation} \label{eq:16ab}
 {}_{0}D_{t}^{1-\varepsilon }y =
 \frac{y(0)}{\Gamma(\varepsilon) t^{1-\varepsilon}}
 +\frac{y^{(1)}(0) t^{\varepsilon}}{\Gamma(1+\varepsilon)}
 +\frac{1}{\Gamma(1+\varepsilon)}
 \int_{0}^{t}(t-\tau)^{\varepsilon}y^{(2)}(\tau)\, d\tau.
 \end{equation}
 From (\ref{eq:16ab}) and (\ref{eq:expansion}) it follows that
 \begin{eqnarray}
 \left.\frac{\partial {}_{0}D_{t}^{\alpha }y}{\partial \alpha}
 \right\vert_{\alpha =1^{-}}
 &\!\!\!\!\!\!\!\!\!=\!\!\!\!& \left.
 -\frac{\partial {}_{0}D_{t}^{1-\varepsilon}y}{\partial \varepsilon}
 \right\vert_{\varepsilon =0^{+}} \notag \\
 &\!\!\!\!\!\!\!\!\!=\!\!\!\!&
 -\frac{y(0)}{t}-y^{(1)}(0) \ln t- \gamma y^{(1)}(t)
 - \int_{0}^{t}y^{(2)}(\tau) \ln(t-\tau)\, d\tau. \label{eq:16bbbb}
 \end{eqnarray}
 Assuming $y(0)=0$ in (\ref{eq:16bbbb}), we recover the results
 presented in \cite{TarasovZaslavsky06} and
 \cite{TofighiGolestani} for the Caputo fractional derivative.
 Since ${}_{0}D_{t}^{\alpha}y={}_{0}^{C}D_{t}^{\alpha}y$
 when $y(0)=0$, it follows that
 with $y(0)=0$ (\ref{eq:16bbbb}) becomes
 $$
 \left.\frac{\partial {}_{0}D_{t}^{\alpha}y}{\partial\alpha}
 \right\vert_{\alpha =1^{-}}
 =
 \left.\frac{\partial\, {}_{0}^{C}D_{t}^{\alpha}y}{\partial\alpha}
 \right\vert_{\alpha =1^{-}}
 = -y^{(1)}(0) \ln t-\gamma y^{(1)}(t)
 -\int_{0}^{t} y^{(2)}(\tau) \ln (t-\tau)\, d\tau.
 $$
 \erem

 \brem \label{rem:lin ops}
 Functional (\ref{eq:fvp-alpha})
 is a special case of the class of functionals with Lagrangians
 depending on linear operators,
 see \cite[p.\ 51]{Filippov}. Indeed, suppose that Lagrangian $L$
 depends on $t$,
 $y$ and $\cL y$, where $\cL:M\to L^p([0,b])$, $p\in[1,+\infty)$,
 is a linear operator acting on a set of admissible
 functions $M$, which is linear, open and dense in $L^p([0,b])$
 (i.e. $\cL$ belongs to $Lin(M,L^p([0,b]))$, the space of continuous,
 linear functions with the uniform norm).
 Suppose that $L$ is continuously differentiable
 with respect to $t$ and $y$ and twice continuously
 differentiable with respect to $\cL y$.
 Moreover, assume that
 function $t\mapsto L(t,y(t), \cL y(t))$, $t\in [0,b]$,
 is continuous, for all $y\in M$.
 Then the Euler-Lagrange equation reads
 $$
 \frac{\pd L}{\pd y} + \cL^\ast \frac{\pd L}{\pd (\cL y)}=0,
 $$
 where $\cL^\ast$ denotes the adjoint operator of $\cL$.
 In case when $\cL$ is the left Riemann-Liouville operator
 ${}_0D_t^\alpha$, with the adjoint ${}_tD_{b}^\alpha$, the
 latter equation coincides with (\ref{eq:ELeqs-y}).

 If instead of $\cL$ one considers a family $\{\cL_\alpha, \alpha\in
 A\}$, where $A=[0,1]$ (or some other interval), and the mapping
 $A\to Lin(M,L^p([0,b]))$, $\alpha\mapsto \cL_\alpha$, is differentiable,
 then a more general problem of finding stationary points with respect
 to $y$ and $\alpha$ can be formulated. In that case,
 one can derive the second stationarity condition similar to
 (\ref{eq:ELeqs-alpha}):
 $$
 \int_0^b \left( \frac{\d L}{\d\cL}\frac{\d \cL}{\d\al} + \frac{\d
 L}{\d\al}\right)dt=0.
 $$
 \erem

 \section{Equivalent problems}
 \label{SecEqProbl}

 In this section we shall give conditions which provide that problems
 (\ref{eq:min}) and (\ref{eq:minmin}) coincide.

 \bprop \label{EqvivPobl}
 Let the Lagrangian $L$ satisfy (\ref{eq:conditions on L}).
 Assume that for every $\alpha\in [0,1]$ there is a unique
 $y^\ast(t,\alpha)\in\cU$, a solution to (\ref{eq:ELeqs-y}), and
 that the mapping $\alpha\mapsto y^\ast(t,\alpha)$ is
 differentiable as a mapping from $[0,1]$ to $\cU$.
 Then
 the problem $\min_{(y,\alpha)\in\cU\times A}I[y,\alpha]$
 is equivalent to
 the problem $\min_{\alpha\in A}(\min_{y\in\cU}I[y,\alpha])$.
 \eprop

 \pr
 As we have shown in Proposition \ref{prop:ELeqs}, any solution to the
 problem $\min_{(y,\alpha)\in\cU\times A}I[y,\alpha]$ satisfies
 system
 (\ref{eq:ELeqs-y})-(\ref{eq:ELeqs-alpha}). It can be
 solved as follows.
 We first solve (\ref{eq:ELeqs-y}) and the corresponding
 boundary conditions to obtain $y^{\ast}=y^{\ast}(t,\alpha)$.
 According to the assumption, the solution $y^\ast$ is unique.
 Then we insert $y^\ast$ in (\ref{eq:ELeqs-alpha}) to obtain
 $\alpha^{\ast}$.
 In this case, functional $I[y,\alpha]$
 becomes a functional depending only on $\alpha$,
 $\alpha\mapsto I[y^\ast(t,\alpha),\alpha]=I[\alpha]$, and
 therefore (\ref{eq:ELeqs-alpha}) transforms to the total
 derivative of $I[\alpha]$ since
 \begin{eqnarray}
 \frac{dI[\alpha]}{d\alpha} &=&
 \left.\frac{dI[\alpha+\varepsilon]}{d\varepsilon}\right\vert_{\varepsilon=0}
 \nonumber\\
 &=& \int_0^1\left[
 \frac{\pd L}{\pd y}\frac{\pd y}{\pd \alpha} +
 \frac{\pd L}{\pd {}_{0}D_{t}^{\alpha}y} \left(
 {}_{0}D_{t}^{\alpha}\left(\frac{\pd y}{\pd\alpha}\right)
 + \left(\frac{\pd}{\pd\alpha}{}_{0}D_{t}^\alpha\right)y
 \right)
 +\frac{\pd L}{\pd\alpha}
 \right]\, dt \nonumber\\
 &=& \int_0^1\left[
 \frac{\pd y}{\pd \alpha} \left(
 \frac{\pd L}{\pd y} + {}_{t}D_{b}^{\alpha}
 \frac{\partial L}{\partial{}_{0}D_{t}^{\alpha}y}\right)
 +
 \frac{\pd L}{\pd {}_{0}D_{t}^{\alpha}y}
 \left(\frac{\pd}{\pd\alpha}{}_{0}D_{t}^\alpha\right)y
 +\frac{\pd L}{\pd\alpha}
 \right]\, dt \nonumber\\
 &=& \int_0^1 \left(
 \frac{\pd L}{\pd {}_{0}D_{t}^{\alpha}y}
 \left(\frac{\pd}{\pd\alpha}{}_{0}D_{t}^\alpha\right)y
 +\frac{\pd L}{\pd\alpha}
 \right)\, dt, \label{eq:oio}
 \end{eqnarray}
 where we used fractional integration by parts formula
 (\ref{eq:frac int by parts}) in the third,
 and equation (\ref{eq:ELeqs-y}) in the last equality.
 This proves the claim.
 \ep

 The following simple assertion is of particular interest:

 \bprop
 Let $L$ satisfy (\ref{eq:conditions on L}).
 Assume that for every $\alpha\in[0,1]$ there exists a unique
 $y_{\alpha}\in\cU$, a solution to the
 fractional variational problem (\ref{eq:fvp-alpha-min}), and that
 $I[y_{\alpha},\alpha]$ is the corresponding minimal value of
 the functional $I$. Assume additionally that
 $$
 \frac{dI}{d\alpha}(y,\alpha)|_{y=y_\alpha}>0,\quad
 \forall y_\alpha\in\cU.
 $$
 Then the minimal, resp.\ maximal
 value of the functional $I[y,\alpha]$ is attained at $\alpha=0$,
 resp.\ at $\alpha=1$.
 \eprop

 \pr
 Under the above assumptions we have that
 $$
 I[y_{0},0]\leq I[y_{\alpha},0]\leq
 I[y_{\alpha},\alpha]\leq I[y_1,\alpha]
 \leq I[y_1,1],
 \quad \forall\alpha\in[0,1],
 $$
 which proves the claim.
 \ep

 \brem
 The same argument can be applied to
 the case when $dI/d\alpha<0$, for any
 fixed $y_\alpha\in\cU$, i.e. when $I$ is an decreasing function
 of $\alpha$, for any fixed $y_\alpha\in\cU$. In that case the
 minimal, resp.\ maximal value of $I$ is at $\alpha=1$, resp.\
 $\alpha=0$.
 \erem

 \section{Examples}\label{SecEx}

 \subsection{Examples with Lagrangians linear in $y$}

 \bex
 Consider the action integral for the
 inertial motion (no force acting) of a material point of the
 form
 \begin{equation} \label{eq:16aa}
 I[y,\alpha]: = \int_{0}^{1} ({}_{0}D_{t}^{\alpha}y)^{2}\,dt,
 \quad (y,\alpha)\in\cU\times A,
 \end{equation}
 where $\cU:=\{y\in \cU_l\,|\, y(0)=0,\, y(1) =1\}$ and $A=[0,1]$.

 Obviously, the minimal value of $I[y,\alpha]$ is zero, and it is
 attained whenever ${}_0D_t^\alpha y=0$. Solutions to
 equation ${}_0D_t^\alpha y=0$ are of the form $y(t)=C\cdot
 t^{1-\alpha}$, $t\in[0,1]$, $C\in\R$ (cf.\ \cite{SamkoKM}).
 All solutions satisfy the Euler-Lagrange equation
 $$
 {}_{t}D_{1}^{\alpha}({}_{0}D_{t}^{\alpha}y)=0.
 $$
 Stationarity condition (\ref{eq:ELeqs-alpha}) reads:
 $$
 \int_{0}^{1} {}_{0}D_{t}^{\alpha}y
 \left(\psi(1-\alpha) {}_{0}D_{t}^{\alpha}y
 -\frac{1}{\Gamma(1-\alpha)}
 \frac{d}{dt} \int_{0}^{t}
 \frac{\ln (t-\tau)y(\tau)}{(t-\tau)^{\alpha}}\,d\tau\right)\,
 dt=0
 $$
 and is automatically satisfied.

 Note that $C\cdot t^{1-\alpha}\in \cU_l$, for all $C\in\R$,
 but only $t^{1-\alpha}\in \cU$. Hence, we conclude that
 $(y^\ast,\alpha^\ast)=(t^{1-\alpha},\alpha)$, $\alpha\in[0,1]$,
 are solutions to the variational problem $I[y,\alpha]\to\min$,
 for $I$ defined by (\ref{eq:16aa}).
 \eex

 \brem
 If $L(t,y(t),{}_{0}D_{t}^{\alpha}y,\alpha)=
 ({}_{0}D_{t}^{\alpha}y)^{2}+(\alpha-\alpha_0)^2$, for a fixed
 $\alpha_0\in(0,1)$, then the problem
 $\int_0^1 L(t,y(t),{}_{0}D_{t}^{\alpha}y,\alpha)\, dt \to\min$
 has a unique minimizer $(y^\ast,\alpha^\ast)=
 (t^{1-\alpha_0},\alpha_0)$.
 \erem

 \bex
 Let the Lagrangian $L$ be of the form
 $$
 L(t,y(t),{}_0D_t^\alpha y,\alpha):= ({}_0D_t^\alpha y)^2- c\cdot y,
 \quad c\in\R,
 $$
 and let $\cU=\{y\in \cU_l\,|\, y(0)=0\}$, $A=[0,1]$,
 for the variational problem:
 $$
 \min_{(y,\alpha)\in\cU\times A}I[y,\alpha]=
 \min_{(y,\alpha)\in\cU\times A} \int_0^1
 (({}_0D_t^\alpha y)^2- c\cdot y(t))\, dt.
 $$
 Equations (\ref{eq:ELeqs-y}) and (\ref{eq:ELeqs-alpha}) become:
 \bea
 {}_tD_1^\alpha {}_0D_t^\al y & = & c, \label{eq:ELyEx2}\\
 \int_0^1 {}_0D_t^\alpha y\cdot \frac{\pd {}_0D_t^\alpha y}{\pd\alpha}\,dt
 & = & 0.\nonumber
 \eea
 Equation (\ref{eq:ELyEx2}) could be solved as follows.
 First, one introduces a substitution $z(t)={}_0D_t^\alpha y$,
 $t\in[0,1]$, and solve ${}_tD_1^\alpha z = c$:
 $$
 z(t)=c \cdot \frac{(1-t)^\al}{\Gamma(1-\al)},
 \quad t\in[0,1],\alpha\in A.
 $$
 Therefore,
 \beq \label{eq:fde}
 {}_0D_t^\al y (t)= c \cdot \frac{(1-t)^\alpha}{\Gamma(1-\al)},
 \quad t\in[0,1],\alpha\in A.
 \eeq
 Recall,
 $$
 {}_0I_t^\alpha y := \frac{1}{\Gamma(\alpha)} \int_0^t
 (t-\tau)^{\alpha} y(\tau)\, d\tau,
 \quad t\in[0,1],\alpha\in A
 $$
 and apply it on the both sides of (\ref{eq:fde}).
 Using \cite[Th.\ 2.4]{SamkoKM}, i.e.
 ${}_0I_t^\alpha ({}_0D_t^\alpha y)= y(t)$,
 one obtains
 \beas
 y(t,\al) &=& \frac{c}{\Gamma(\al)\Gamma(1+\al)} \int_0^t (t-\tau)^{\al-1}
 (1-\tau)^{\al}\,d\tau\\
 &=& \frac{c}{\Gamma (1+\al)}\sum_{p=0}^{\infty}
 \frac{\Gamma(p-\al)\Gamma(1+p)}{\Gamma(-\al)p\,!\Gamma
 (1+p+\al)}t^{p+\al},
 \quad t\in[0,1],\alpha\in A.
 \eeas
 This solution is unique and belongs to $\cU$.
 Since $\al\mapsto y(t,\al)$ is differentiable,
 Proposition \ref{EqvivPobl} holds.

 We substitute obtained $y(t,\al)$ into $I[y,\al]$ which yields
 \begin{eqnarray*}
 I[\al] &=&
 \int_0^1 \left(\left(c \cdot \frac{(1-t)^\al}{\Gamma(1-\al)}\right)^2-
 \frac{c^2}{\Gamma(\al)\Gamma(1+\al)}\int_0^t(t-\tau)^{\al-1}
 (1-\tau)^{\al}d\tau\right)\,dt \\
 &=&
 \int_0^1 \left(\left(c \cdot \frac{(1-t)^\al}{\Gamma(1-\al)}\right))^2 -
 \frac{c^2}{\Gamma (1+\al)}\sum_{p=0}^{\infty}
 \frac{\Gamma(p-\al)\Gamma(1+p)}{\Gamma(-\al)p\,!\Gamma
 (1+p+\al)}t^{p+\al}\right)\,dt.
 \end{eqnarray*}
 Simple numerical calculations show that $I[\alpha]$ is an increasing
 function and attains extremal values at the boundaries.

 \brem
 Equation (\ref{eq:ELyEx2}) represents a fractional generalization of
 the equation of motion for a material point (with unit mass)
 under the action of constant force equal to $c$. Our result
 shows that an optimal value of Hamilton's action is attained for
 $\al=1$, that is for integer order dynamics. We note that
 different generalizations of classical equation of motion can be
 found in \cite{Kwok}, where the problem
 ${}_0D_t^\al y = c$, $1<\alpha\leq 2$, was analyzed.
 \erem
 \eex

 \subsection{Examples with Lagrangians linear in ${}_0D_t^\alpha y$}

 \bex
 Let
 \beq \label{eq:ex3-L}
 L(t,y(t),{}_0D_t^\alpha y,\alpha):=
 \Gamma(1-\alpha){}_0D_t^\alpha y - \frac{1}{2} c y^2,
 \quad c>0, c\not=1,
 \eeq
 and consider the problem of finding stationary points for the
 functional (\ref{eq:fvp-alpha}),
 where $\cU:=\{y\in\cU_l\,|\,y(0)=\frac{1}{c}\}$
 and $\alpha_0<\frac12$.
 Note that $L$ satisfies the so called primary constraint in Dirac's
 classification of systems with constraints
 (cf.\ \cite{HenneauxTeitelboim}).
 In the setting of fractional derivatives such
 Lagrangians have been recently treated in \cite{Bal04} and
 \cite{MuslihBaleanu}.

 Equations (\ref{eq:ELeqs-y}) and (\ref{eq:ELeqs-alpha}) become
 \beq \label{eq:ex3-ELy}
 \Gamma(1-\alpha){}_tD_1^\alpha 1 - c y=0
 \eeq
 and
 \beq \label{eq:ex3-ELalph}
 \int_0^1 \left(\Gamma(1-\alpha)\frac{\pd{}_0D_t^\alpha y}{\pd\alpha}
 + \frac{\pd \Gamma(1-\alpha)}{\pd\alpha}\right)\, dt=0.
 \eeq
 Equation (\ref{eq:ex3-ELy}) has a unique solution
 $y^\ast=\frac{1}{c(1-t)^\alpha}$, $t\in[0,1]$,
 $\alpha\in[0,\alpha_0]$. This implies
 \beas
 I[y^\ast,\alpha] &=& \int_0^1 \left[\frac{d}{dt} \int_0^t
 \frac{d\tau}{c(1-\tau)^\alpha(t-\tau)^\alpha}
 -\frac{1}{2c(1-t)^{2\alpha}}\right] \, dt \\
 &=& \left. \int_0^t
 \frac{d\tau}{c(1-\tau)^\alpha(t-\tau)^\alpha}
 \right\vert_{0}^{1}
 - \int_0^1
 \frac{1}{2c(1-t)^{2\alpha}} \, dt \\
 &=& \int_0^1
 \frac{d\tau}{c(1-\tau)^{2\alpha}}
 - \int_0^1
 \frac{1}{2c(1-t)^{2\alpha}} \, dt \\
 &=&
 \frac{1}{2c} \int_0^1
 \frac{1}{(1-t)^{2\alpha}} \, dt.
 \eeas
 Since $\alpha_0<1/2$ we have that $I[y^\ast,\alpha]$ exists and
 is an increasing function with respect to $\alpha$. Hence,
 $I[y^\ast,\alpha]$ attains its
 minimal value at $\alpha=0$, and it equals $1/(2c)$.
 We also have that the maximal value of $I[y^\ast,\alpha]$ is
 attained at $\alpha_0$.
 \eex

 \bex
 Let $\cU:=\{y\in \cU_l\,|\, y(0)=0\}$,
 $c\neq 0$ and let $L$ be of the form
 \begin{equation}\label{eq:exL}
 L(t,y(t),{}_0D_t^\alpha y,\alpha): =
 c\cdot{}_0D_t^\alpha y + f(y(t)),\quad t\in[0,1],
 \end{equation}
 where the properties of $f$ are going to be specified.

 In this example we are dealing with integrable functions which
 can take values $+\infty$ or $-\infty$ at some points.
 We are going to analyze stationary points of
 $$
 I[y,\alpha]= \int_0^1
 (c\cdot{}_0D_t^\alpha y(t) + f(y(t)))\, dt.
 $$

 Equations (\ref{eq:ELeqs-y}) and (\ref{eq:ELeqs-alpha}) become
 \bea
 {}_tD_1^\alpha c + \frac{\pd f}{\pd y} & = & 0 \label{eq:ELyEx}\\
 c \cdot \int_0^1 \frac{\pd {}_0D_t^\alpha y}{\pd\alpha}\,dt
 & = & 0. \label{eq:ELalEx}
 \eea
 Since ${}_tD_1^\alpha c =
 \frac{c}{\Gamma (1-\alpha) (1-t)^{\alpha}}$, $t\in[0,1]$,
 we see that in order to solve (\ref{eq:ELyEx})-(\ref{eq:ELalEx})
 we have to assume that $f\in\cC^1(\R)$, and that $f'$ is
 invertible so that $t\mapsto (f')^{-1}
 (\frac{c}{\Gamma (1-\alpha) (1-t)^{\alpha}})\in\cU_l$.
 Then equation (\ref{eq:ELyEx}) is solvable
 with  respect to $y$:
 \begin{equation}\label{eq:Ex_y_od_c}
 y_c(t,\alpha) = \left(\frac{\d f}{\d y}\right)^{-1}
 \left( \frac{c}{\Gamma (1-\alpha) (1-t)^{\alpha}} \right), \quad t\in [0,1].
 \end{equation}
 Since $c\neq 0$, (\ref{eq:ELalEx}) implies
 $\int_0^1 \frac{\pd {}_0D_t^\alpha y}{\pd\alpha}\,dt=0$. Thus
 \bea
 0 & = &
 \int_0^1 \frac{\pd {}_0D_t^\alpha y}{\pd\alpha}\,dt
 = \int_0^1 G(y, \alpha)(t)\,dt
 = \int_0^1 \frac{d}{dt}(f_1\ast_t y)(t,\alpha)\,dt \nonumber \\
 & = &  f_1\ast_t y(t,\alpha)|_{t=1} - f_1\ast_t y(t,\alpha)|_{t=0} =
 f_1\ast_t y(t,\alpha)|_{t=1},\label{integral}
 \eea
 where we have used that  $f_1 \in L^1([0,1])$ and that $y\in \cU$.
 Substitution of (\ref{eq:Ex_y_od_c}) into (\ref{integral}) gives
 $(f_1\ast y_c(t,\alpha))(t)|_{t=1}=0$ or
 \begin{equation} \label{eq:exc}
 \int_0^1 \frac{\psi(1-\alpha)-\ln (1-\tau)}{\Gamma (1-\alpha)(1-\tau)^{\alpha}}
 \left(\frac{\d f}{\d y}\right)^{-1}\left(\frac{c}{\Gamma (1-\alpha)(\tau-1)^{\alpha}}\right)
 \,d\tau=0.
 \end{equation}
 Solving this equation is obviously difficult.
 Hence, we consider some special cases.

 a)\ $f(y(t)):= d \cdot \frac{y(t)^2}{2}$, $t\in[0,1]$, $d\in\R$.
 Then the Lagrangian is
 $$
 L(t,y(t),{}_0D_t^\alpha y,\alpha) =
 c\cdot{}_0D_t^\alpha y(t) + d \cdot \frac{y(t)^2}{2},
 $$
 and
 $$
 y_c(t,\alpha)=-\frac{1}{d}\frac{c}{\Gamma(1-\al)(1-t)^\al},
 \quad t\in[0,1].
 $$
 Also, (\ref{eq:exc}) becomes
 $$
 \int_0^1 \frac{\psi(1-\alpha)-\ln (1-\tau)}
 {\Gamma (1-\alpha)^2(1-\tau)^{2\alpha}}\,d\tau =0.
 $$
 By a simple numerical calculation one shows that this equation
 does not have any solution for $\al\in (0,1)$. Hence, in this
 case there does not exist any point $(y,\alpha)$ which is an
 extremal of functional $I[y,\alpha]$.

 b)\ $f(y(t)): = \ln y(t)$, $t\in[0,1]$. Then (\ref{eq:ELyEx})
 becomes
 $$
 \frac{c}{\Gamma(1-\al)(1-t)^\al}=\frac{1}{y},
 $$
 and therefore
 $$
 y=\frac{\Gamma(1-\al)}{c}(1-t)^\al\in AC([0,1]).
 $$
 In this particular case we take the set of admissible functions
 to be $\cU:=\{y\in\cU_l\,|\,y(1)=0\}$.
 Using (\ref{integral}), equation (\ref{eq:ELalEx}) reads
 $$
 \int_0^1 \left(\psi (1-\al)-\ln (1-\tau)\right) \,d\tau = 0,
 $$
 which, after integration, yields $\psi (\al-1)=1$.
 A unique solution of this equation in $(0,1)$ is
 $\alpha=0.604...$.

 Therefore, a unique stationary point of
 $$
 I[y,\alpha]=\int_0^1 (c\cdot{}_0D_t^\alpha y
 + \ln y(t))\, dt
 $$
 is the point $(y,\alpha)=
 (\frac{\Gamma(0,396)}{c}(1-t)^\al;0,604)$.
 \eex

 \brem
 So far, we have considered variational problems defined via functionals
 of type (\ref{eq:fvp-alpha}). In fact, we have allowed fractional
 derivatives of functions to appear in Lagrangians. The natural
 generalization of such problems consists of replacing the
 Lebesgue integral in (\ref{eq:fvp-alpha}) by the Riemann-Liouville
 fractional integral.
 More precisely, for $\beta>0$ set
 \beas
 I_\beta[y,\alpha] :&=& {}_0I_b^\beta
 L(t,y(t),{}_{0}D_{t}^{\alpha}y,\alpha) \\
 &=& \frac{1}{\Gamma (\beta)}
 \int_{0}^{b}
 (b-t)^{\beta -1}L(t,y( t),{}_{0}D_{t}^{\alpha}y,\alpha)\, dt,
 \quad t\in (0,b).
 \eeas
 Then the fractional variational problem consists of finding extremal
 values of the functional $I_\beta[y,\alpha]$.

 In the above construction we have used the left Riemann-Liouville
 fractional integral of order $\beta$ (which, in general, differs
 from the order of fractional differentiation $\alpha$), evaluated
 at $t=b$. The choice $\beta=1$ turns us back to the problem
 (\ref{eq:fvp-alpha}).

 The study of such fractional variational problems is reduced
 to the case we have already considered in the following way.
 It suffices to redefine the Lagrangian as
 $$
 L_1(t,y(t),{}_{0}D_{t}^{\alpha}y,\alpha,\beta):=
 \frac{1}{\Gamma (\beta)} (b-t)^{\beta -1}
 L(t,y( t),{}_{0}D_{t}^{\alpha}y,\alpha).
 $$
 Then we have to consider the functional
 \beq \label{eq:fracfrac-vp}
 I_\beta[y,\alpha]=\int_{0}^{b}
 L_1(t,y(t),{}_{0}D_{t}^{\alpha}y,\alpha,\beta)\, dt.
 \eeq
 In case $\beta>1$, $L_1$ is of the same regularity as $L$, so
 the straightforward application of the results derived in
 previous sections to the Lagrangian $L_1$ leads to
 the optimality conditions for the variational problem
 defined via the functional (\ref{eq:fracfrac-vp}).
 However, when $0<\beta<1$, continuity, as well as
 differentiability of $L_1$ with respect to $t$ may be violated
 (which depends of the explicit form of $L$), and hence it may be
 not possible to use the theory developed so far.
 \erem

 \section{Conclusion}

 We formulated Hamilton's principle so that the order of derivative
 in the Lagrangian is also subject to variation. The stationarity
 conditions are derived in (\ref{eq:ELeqs-y}) and
 (\ref{eq:ELeqs-alpha}). We introduced additional assumptions which
 resulted in equivalent problems, simpler for solving.
 Several examples are given in order to illustrate the theory
 presented in the paper.
 We concluded our work with a consideration
 of Hamilton's principle defined in terms of Riemann-Liouville
 fractional integrals.

 \subsection*{Acknowledgement}

 This work is supported by Projects $144016$ and $144019$ of the
 Serbian Ministry of Science and START-project Y-237 of the
 Austrian Science Fund.



 \end{document}